\documentclass{amsart}

\usepackage{amssymb}
\usepackage{mathtools}
\usepackage[pdftex,hypertexnames=false]{hyperref}
\usepackage{tikz-cd}
\usepackage[shortlabels]{enumitem}

\newtheorem{theorem}{Theorem}[section]
\newtheorem{lemma}[theorem]{Lemma}
\newtheorem{proposition}[theorem]{Proposition}
\newtheorem{corollary}[theorem]{Corollary}

\theoremstyle{definition}
\newtheorem{definition}[theorem]{Definition}
\newtheorem{example}[theorem]{Example}

\theoremstyle{remark}
\newtheorem{remark}[theorem]{Remark}

\numberwithin{equation}{section}

\newcommand{\cat}[1]{\mathsf{#1}}

\newcommand{\CompHaus}{\cat{CH}}

\newcommand{\CompMet}{\cat{CompMet}}
\newcommand{\OKH}{\cat{OrdCH}}
\newcommand{\MKH}{\cat{MetCH}}

\newcommand{\Met}{\cat{Met}}

\newcommand{\BRight}[2]{B^{\vartriangleright}(#1,#2)}
\newcommand{\BLeft}[2]{B^{\vartriangleleft}(#1,#2)}
\newcommand{\BRightC}[2]{\overline{B}^{\vartriangleright}(#1,#2)}
\newcommand{\BLeftC}[2]{\overline{B}^{\vartriangleleft}(#1,#2)}
\newcommand{\BBoth}[2]{B(#1,#2)}

\DeclareMathOperator{\ev}{ev}

\DeclareMathOperator{\upc}{\uparrow\!}
\DeclareMathOperator{\downc}{\downarrow\!}

\newcommand{\doo}[1]{\overset{\centerdot}{#1}}

\author[D. Hofmann]{Dirk Hofmann}
\address{CIDMA, Department of Mathematics, University of Aveiro, 3810-193
Aveiro, Portugal}
\email{dirk@ua.pt}

\author[R. Prezado]{Rui Prezado}
\address{CMUC, Department of Mathematics, University of Coimbra, 3000-143
Coimbra, Portugal}
\email{ruiprezado@mat.uc.pt}



\title[Exponentiable Metric Compact Hausdorff spaces]{Compact metric spaces as
exponentiable metric compact Hausdorff spaces}

\begin{document}

\begin{abstract}
  The class of compact metric spaces (metric spaces which induce a compact
  topology), despite being an important class of metric spaces, is not closed
  under many natural constructions. For instance, infinite metric products of
  compact metric spaces will often fail to induce a compact topology, while on
  the other hand, the topological product of any family of compact metric
  spaces is always a compact topological space, evidencing a discrepancy
  between topological and metric products.

  The notion of metric compact Hausdorff space -- compact Hausdorff spaces
  equipped with a lower semi-continuous metric -- provides a class of spaces
  which keep track of metric data and topological data separately, thereby
  forming a complete and cocomplete category, while including all compact
  metric spaces.

  Our main contribution is that this inclusion can be characterised via
  universal property: the compact metric spaces are the exponentiable metric
  compact Hausdorff spaces with respect to a natural monoidal structure on the
  respective category.
\end{abstract}

\maketitle

\section*{Introduction}
\label{rui:sect:intro}
A fruitful theory of (topological) spaces should allow for constructing new
spaces out of a family of old ones; it is also often desirable for two given
constructions to commute with one other (see also \cite{BHL91}).  Prominent
examples of this phenomena are the so-called ``convenient categories of
spaces'' \cite{Ste67}, that is, categories of spaces which are Cartesian
closed (in the sense that they do admit a canonical function space \([X,Y]\)
for all spaces \(X\) and \(Y\)), for this property implies that ``products
commute with all quotients''.

As it is well-known, many important categories of topology are \emph{not}
Cartesian closed, motivating the study of those spaces \(X\) -- called
\emph{exponentiable} spaces -- which do admit such a function space \([X,Y]\)
for all \(Y\); put more formally, those spaces \(X\) such that the endofunctor
\(-\times X\) has a right adjoint \([X,-]\).  For instance, the exponentiable
topological spaces are precisely the core-compact spaces \cite{DK70}, and the
exponentiable metric spaces are the totally quasi-convex spaces \cite{CH06}.
In this context, it might seem rather disappointing that the exponentiable
compact Hausdorff spaces are trivial: as shown in \cite{CM91}, a compact
Hausdorff space is exponentiable in the category \(\CompHaus\) of compact
Hausdorff spaces and continuous maps if and only if it is finite.  The
situation does not improve if we move to the category \(\OKH\) of Nachbin's
ordered compact Hausdorff spaces with monotone and continuous maps
\cite{Nac65}: it is not difficult to prove that in this category the
exponentiable objects are essentially the finite ones as well (see
Corollary~\ref{d:cor:1}).

Surprisingly, we have more luck considering \emph{metric compact Hausdorff
spaces} and appropriate morphisms.  These spaces constitute a simultaneous
generalisation of ordered compact Hausdorff spaces and of compact metric
spaces; in a nutshell, a metric compact Hausdorff space is given by a set
equipped with a metric and a compact Hausdorff topology (not necessarily
induced by the metric) which are compatible in an appropriate sense.  We
stress that every compact metric space is also a metric compact Hausdorff
space, and several arguments usually applied to compact metric space can be
used in this context as well (see, for instance, Proposition~\ref{d:prop:3}).
This fact makes us confident that metric compact Hausdorff spaces provide a
flexible generalisation of compact metric spaces.  Metric compact Hausdorff
spaces and the corresponding category were introduced in \cite{HR18} based on
\cite{Tho09}. However, with different motivation, similar combinations of
topology and metric were also investigated in \cite{Mat00, BYaa08, BYaa13,
PST26}.

Our main contribution is that the inclusion of compact metric spaces into
metric compact Hausdorff spaces can be characterised via universal property:
the compact metric spaces are the exponentiable metric compact Hausdorff space
with respect to a natural monoidal structure on the category of metric compact
Hausdorff spaces and morphisms.  We start by characterising those metric
spaces which are exponentiable with respect to a suitable monoidal structure
on the category of metric spaces.  We obtain an easily checkable
characterization (Theorem~\ref{thm:blob-convexity}) which allows us to prove
that exponentiability for certain Minkowski tensors is equivalent to total
quasi-convexity (Corollary~\ref{cor:minkowski}).  Building on this, we obtain
our main result (Theorem \ref{d:thm:4}): the exponentiable (with respect to a
suitable monoidal structure) separated metric compact Hausdorff spaces are
precisely the compact metric spaces.  Hence, upholding the importance of
metric compact Hausdorff spaces, the main result of our paper shows that
compact metric spaces arise naturally as an important subclass of all
separated metric compact Hausdorff spaces.

\subsection*{Funding}

The authors acknowledge the partial support by CIDMA under the Portuguese
Foundation for Science and Technology (FCT, \url{https://ror.org/00snfqn58})
Multi-Annual Financing Program for R\&D Units, grants \url{UID/4106/2025} and
\url{UID/PRR/4106/2025}.  The second author acknowledges the partial support
by the Center for Mathematics of the University of Coimbra (CMUC,
\url{https://doi.org/10.54499/UID/00324/2025}) under the Portuguese Foundation
for Science and Technology (FCT), Grants \url{UID/00324/2025} and
\url{UID/PRR/00324/2025}.

\section{Preliminaries}
In this paper we use the designation \emph{metric} in the generalised sense of
\cite{Law73}; namely as a distance function
\begin{displaymath}
  d \colon X\times X \longrightarrow [0,\infty]
\end{displaymath}
such that \(d(x,x)=0\) and \(d(x,y)+d(y,z)\geq d(x,z)\), for all \(x,y,z\in
X\).  A \emph{metric space} is a pair \((X,d)\) consisting of a set \(X\) and
a metric \(d\) on \(X\), but we often write simply \(X\) instead of \((X,d)\).
A \emph{metric map} \(f\colon X\to Y\) between metric spaces consists of a
function \(f \colon X \to Y \) between the sets \(X\) and \(Y\) satisfying
\(d(x_1, x_2) \geq d(f(x_1), f(x_2))\), for all \(x_1, x_2\in X\).  The
category of metric spaces and metric maps is denoted by \(\Met\).

\begin{example}
  \label{d:ex:1}
  Of particular importance is the metric space \([0,\infty]\) with the metric
  \(\mu\) defined by \(\mu(u,v)=\max\{v-u,0\}\), for all \(u,v\in
  [0,\infty]\).
\end{example}

Since our notion of metric is \emph{asymmetric}; that is, we generally have
\( d(x,y) \neq d(y,x) \) for some pairs \( x, y \) in a metric space \( X \),
we have two distinct notions of open balls and closed balls, denoted by
\begin{align*}
  \BLeft{x}{r} &= \{ y \mid d(y,x) < r \}
  & \BRight{x}{r} &= \{ y \mid d(x,y) < r \} \\
  \BLeftC{x}{r} &= \{ y \mid d(y,x) \leq r \}
  & \BRightC{x}{r} &= \{ y \mid d(x,y) \leq r \}
\end{align*}
which are respectively called the \emph{lower} and \emph{upper} open/closed
balls with center \( x \in X \) and radius \( r \in [0,\infty] \).  By
definition, the lower and upper open balls form a sub-basis for the
\emph{topology on \(X\) induced by \(d\)}, hence, the collection of all
\emph{symmetric open balls}
\begin{displaymath}
  \BBoth{x}{r}=\BLeft{x}{r}\cap\BRight{x}{r}=\{y\mid d(x,y)<r \:\&\: d(y,x)<r\}
\end{displaymath}
with center \(x\in X\) and radius \(r\in [0,\infty]\) form a basis for this
topology.  As usual, a metric space is called \emph{compact} whenever this
topology is compact.  We also note that this topology is Hausdorff if and only
if \(X\) is \emph{separated}, that is,
\begin{displaymath}
  (0=d(x,y)\;\&\; 0=d(y,x)) \implies x=y
\end{displaymath}
for all \(x,y\in X\).

The category \(\Met\) has a terminal object, namely the one-element space
\(1\), and for metric spaces \( X, Y \), the product \(X\times Y\) of \(X\)
and \(Y\) is given by the Cartesian product \( X \times Y \) of the underlying
sets equipped with the \emph{max-metric} defined by
\begin{displaymath}
  d((x_1,y_1),(x_2,y_2))=\max\{d(x_1,x_2),d(y_1,y_2)\},
\end{displaymath}
for all \(x_1,x_2\in X\) and \(y_1,y_2\in Y\).

A metric space \( X \) is \emph{exponentiable} \emph{in \( \Met \)} whenever
the functor \( - \times X \colon \Met \to \Met \) has a right adjoint, which
we denote by \( (-)^X \).  The category \(\Met\) is not Cartesian closed, and
below (Theorem~\ref{d:thm:2}) we recall the characterisation of exponentiable
metric spaces obtained in \cite{CH06} (see also \cite{CHS09}) via a
\textit{quasi-convexity} condition:

\begin{definition}
  Let \( X \) be a metric space.
  \begin{enumerate}
    \item
      We say \( X \) is \emph{totally quasi-convex} whenever, for all \( x_1,
      x_2 \in X \) and all \( u_1, u_2 \in [0,\infty] \) such that \(
      d(x_1,x_2) = u_1+u_2 \), for all \( \varepsilon > 0 \), there exists \(
      z \in X \) such that \( d(x_1,z) \leq u_1+\varepsilon \) and \( d(z,x_2)
      \leq u_2+\varepsilon \).
    \item
      We say \( X \) is \emph{totally convex} if for all \( x, y \in X \) and
      all \(u_1.u_2 \in [0,\infty] \) with \( d(x,y) = u_1+u_2 \), there
      exists \( z \in X \) such that \( d(x,z) \leq u_1 \) and \( d(z,y) \leq
      u_2 \).
  \end{enumerate}
\end{definition}

\begin{remark}
    \label{d:rem:1}
    Note that, thanks to the triangle inequality, one necessarily has
    \(d(x,z)=u_1\) and \(d(z,y)=u_2\) when \(u_1+u_2<\infty\) in the
    definition of a totally convex space.  Clearly, every totally convex
    metric space is totally quasi-convex, and one easily verifies the converse
    for compact metric spaces.\footnote{In \cite{CH06} we stated erroneously
    that every complete totally quasi-convex metric space is totally convex, a
    mistake brought to our attention by Carlyle Stewart (University of
    Stellenbosch, South Africa).}
\end{remark}

\begin{theorem}[\cite{CH06}]
  \label{d:thm:2}
  The following assertions are equivalent, for a metric space \(X\).
  \begin{enumerate}[(i)]
    \item
      \(X\) is exponentiable in \(\Met\).
    \item
      For all \( x_1, x_2 \in X \) and all \( u_1, u_2 \in [0,\infty]\),
      \begin{displaymath}
        \inf_{z\in X}\big(\max\{d(x_1,z),u_1\}+\max\{d(z,x_2),u_2\}\big)
        \leq \max\{d(x_1,x_2),u_1+u_2\}.
      \end{displaymath}
    \item
      \( X \) is totally quasi-convex.
  \end{enumerate}
\end{theorem}

For metric spaces \(X\) and \(Y\), another important metric on the Cartesian
product \(X\times Y\) of the sets \(X\) and \(Y\) is the \emph{plus-metric}
defined by
\begin{displaymath}
  d((x_1,y_1),(x_2,y_2))\coloneq d(x_1,x_2)+d(y_1,y_2),
\end{displaymath}
for all \(x_1,x_2\in X\) and \(y_1,y_2\in Y\).  We denote the resulting metric
space by \(X\oplus Y\).  In contrast to the categorical product, this
construction makes \(\Met\) a closed symmetric monoidal category (see
\cite{Law73}); here, the right adjoint \([X,-]_+\) of \(-\oplus X\) is given
by
\begin{displaymath}
  [X,Y]_+  = \{f \colon X\to Y\mid f\text{ is metric}\},\quad
  d(f,g)  =\sup_{x\in X}d(f(x),g(x))
\end{displaymath}
for all metric spaces \(X\) and \(Y\) and metric maps \(f,g \colon X\to Y\).

Below we investigate the interaction of exponentiability with two well-known
constructions of metric spaces.  Here we will use here the following notation:
for a metric space \(X\), \(A,B\subseteq X\) and \(\varepsilon>0\),
\begin{displaymath}
  C_{A,B,\varepsilon} = \{z\in X\mid \exists x\in A, y\in B\,.\,
      (d(x,z)\leq u_1+\varepsilon\;\&\;d(z,y)\leq u_2+\varepsilon)\}.
\end{displaymath}

First, we consider the metric space \(\mathcal{H}(X)\) of all subsets of
\(X\), with the \emph{asymmetric Hausdorff metric} defined by
\begin{displaymath}
  d(A,B)=\sup_{y\in B}\inf_{x\in A}d(x,y),
\end{displaymath}
for all \(A,B\subseteq X\) (\cite{Pom05, Hau14}, but see also \cite{LLP93}.
We note that, since \(d(x,y)=d(\{x\},\{y\})\) for all \(x,y\in X\), the
canonical map \(X\to \mathcal{H}(X),x\mapsto\{x\}\) is an isometric embedding.

\begin{proposition}
  Let \(X\) be a metric space.  Then \(X\) is totally quasi-convex if and only
  if \(\mathcal{H}(X)\) is totally quasi-convex.
\end{proposition}
\begin{proof}
  Assume first that \(X\) is totally quasi-convex.  Let \(A,B\subseteq X\) and
  \(u_1,u_2\in [0,\infty]\) with \(d(A,B)=u_1+u_2\).  Let \(\varepsilon>0\).
  One easily verifies that
  \begin{displaymath}
    d(A,C_{A,B,\varepsilon})\leq u_1+\varepsilon
    \quad\text{and}\quad
    d(C_{A,B,\varepsilon},B)\leq u_2+\varepsilon.
  \end{displaymath}
  Assume now that \(\mathcal{H}(X)\) is totally quasi-convex.  Let
  \(x_1,x_2\in X\) and \(u_1,u_2\in [0,\infty]\) with \(d(x_1,x_2)=u_1+u_2\).
  Let \(\varepsilon>0\).  By hypothesis, there exists a subset \(C\subseteq
  X\) so that
  \begin{displaymath}
    u_1+\varepsilon>u_1+\frac{\varepsilon}{2}\geq d(\{x_1\},C)=\sup_{z\in C}d(x_1,z)
  \end{displaymath}
  and
  \begin{displaymath}
    u_2+\varepsilon>u_2+\frac{\varepsilon}{2}\geq d(C,\{x_2\})=\inf_{z\in C}d(z,x_2).
  \end{displaymath}
  By the second inequality, there is some \(z_0\in C\) with \(d(z_0,x_2)\leq
  u_2+\varepsilon\), by the first inequality, we also have \(d(x_1,z_0)\leq
  u_1+\varepsilon\).
\end{proof}

For a metric space \(X\), the set \(\mathcal{U}(X)\) of all ultrafilters on
\(X\) becomes a metric space with distance
\begin{displaymath}
  d(\mathfrak{x}, \mathfrak{y})
    = \sup_{A\in \mathfrak{x},B\in \mathfrak{y}} \inf_{x\in A, y\in B} d(x,y),
\end{displaymath}
for all \(\mathfrak{x}, \mathfrak{y}\in \mathcal{U}(X)\) (see
\cite[Section~III.5]{HST14}). In particular, for all \(x,y\in X\), we have
\(d(\doo{x},\doo{y})=d(x,y)\), making the map \(X \to \mathcal{U}(X),\,x
\mapsto\doo{x}\) an isometric embedding. We also note that the metric space
\(\mathcal{U}(X)\) is in general not compact: for instance, if \(X\) is
discrete, then \(\mathcal{U}(X)\) is discrete.

\begin{proposition}
  \label{d:prop:2}
  Let \(X\) be a metric space. Then the following assertions are equivalent.
  \begin{enumerate}[(i)]
    \item
      \(X\) is exponentiable, i.e.\ totally quasi-convex.
    \item
      \(\mathcal{U}(X)\) is totally convex.
    \item
      \(\mathcal{U}(X)\) is totally quasi-convex.
  \end{enumerate}
\end{proposition}
\begin{proof}
  Assume first that \(X\) is exponentiable.  Let \(\mathfrak{x}, \mathfrak{y}
  \in \mathcal{U}(X)\) with \(d(\mathfrak{x}, \mathfrak{y})=u_1+u_2\).  The
  collection
  \begin{displaymath}
    \{C_{A,B,\varepsilon}\mid
    A\in \mathfrak{x}, B\in \mathfrak{y}, \varepsilon>0\}
  \end{displaymath}
  is a filter base on \(X\), which can be extended to an ultrafilter
  \(\mathfrak{z}\in \mathcal{U}(X)\).  For each \(A\in \mathfrak{x}\), \(C\in
  \mathfrak{z} \) and \( \varepsilon>0 \), we have \(C\cap C_{A,X,\varepsilon}
  \neq \varnothing\), hence, there are \(x\in A\) and \(z\in C\) with
  \(d(x,z)\leq u_1+\varepsilon\).  This proves \(d(\mathfrak{x}, \mathfrak{z})
  \leq u_1 + \varepsilon\), for each \(\varepsilon>0\), and therefore
  \(d(\mathfrak{x}, \mathfrak{z})\leq u_1\).  A similar argument shows that
  \(d(\mathfrak{z}, \mathfrak{y})\leq u_2\). Therefore \(\mathcal{U}(X)\) is
  totally convex.

  Clearly, if \(\mathcal{U}(X)\) is totally convex, then \(\mathcal{U}(X)\) is
  totally quasi-convex.

  Assume now that \(\mathcal{U}(X)\) is totally quasi-convex and let \(x,y\in
  X\) with \(d(x,y)=u_1+u_2\) and let \(\varepsilon>0\).  By hypothesis, there
  exists \(\mathfrak{z}\in \mathcal{U}(X)\) such that \(d(\doo{x},
  \mathfrak{z})\leq u_1+\frac{\varepsilon}{2}\) and
  \(d(\mathfrak{z},\doo{y})\leq u_2+\frac{\varepsilon}{2}\).  Therefore, for
  every \(C\in \mathfrak{z}\),
  \begin{displaymath}
    \inf_{z\in C}d(x,z)<u_1+\varepsilon
    \quad\text{and}\quad
    \inf_{z\in Z}d(z,y)<u_2+\varepsilon.
  \end{displaymath}
  Put \(D=\{z\in Z\mid d(x,z)\leq u_1+\varepsilon\}\).
  The first inequality above implies that \(D\in \mathfrak{z}\), and then the
  second one entails that there exits some \(z\in D\) with \(d(z,y)\leq
  u_2+\varepsilon\). We conclude that \(X\) is totally quasi-convex.
\end{proof}

In the remainder of this section we will consider the main objects of interest
for this paper: metric spaces equipped with an additional compact Hausdorff
topology.

\begin{definition}
  A \emph{metric compact Hausdorff space} \(X\) consists of a compact
  Hausdorff space \(X\) together with a lower semicontinuous metric
  \begin{displaymath}
    d \colon X\times X \longrightarrow [0,\infty],
  \end{displaymath}
  that is, the metric \(d\) is continuous with respect to the topology on
  \([0,\infty]\) generated by the basic closed sets \([0,u]\) (\(u\in
  [0,\infty]\)).
\end{definition}

\begin{example}
  The metric space \([0,\infty]\) of Example~\ref{d:ex:1} becomes a metric
  compact Hausdorff space with the \emph{Euclidean topology} on
  \([0,\infty]\), that is, the topology generated by the basic closed sets
  \([0,u]\) and \([u,\infty]\) (\(u\in [0,\infty]\)).  We note that this
  topology differs from the one induced by the metric \(\mu\), which
  is not compact.
\end{example}

\begin{example}
  For each metric space \(X\), the metric space \(\mathcal{U}(X)\) becomes a
  metric compact Hausdorff space with the \emph{Zariski topology} (also called
  \emph{hull-kernel topology}, see \cite{GJ60}), that is, the topology with
  basic open sets
  \begin{displaymath}
    \{\mathfrak{x}\in \mathcal{U}(X)\mid A\in \mathfrak{x}\}\quad (A\subseteq X).
  \end{displaymath}
  We find it worthwhile to point out that this topology only depends on the
  underlying set of \(X\) and not on its metric.
\end{example}

\begin{example}
  \label{d:ex:2}
  Recall from \cite{Nac65} that an \emph{ordered compact Hausdorff space} is a
  compact Hausdorff space \(X\) equipped with an order relation \(\leq\)
  (i.e.\ a reflexive and transitive relation on \(X\)) such that
  \begin{displaymath}
    \{(x,x')\mid x\leq x'\}\subseteq X\times X
  \end{displaymath}
  is closed in the product space \(X\times X\).

  Every ordered compact Hausdorff space can be viewed as a metric compact
  Hausdorff space via the metric \(d\) defined by
  \begin{displaymath}
    d(x,y)=
    \begin{cases}
      0 & \text{if \(x\leq y\)},\\
      \infty & \text{if \(x\nleq y\)},
    \end{cases}
  \end{displaymath}
  for all \(x,y\in X\).
\end{example}

\begin{remark}
  \label{d:rem:2}
  If \(X\) be a metric compact Hausdorff space, then the closed balls are
  closed but the open balls need not be open (see Proposition~\ref{d:prop:1}).
  Furthermore, if \(X\) is separated, then the metric topology refines the
  topology of \(X\) (see \cite[Lemma~1.4]{BYaa08}).
\end{remark}

The following result is a slight generalisation of Remark~\ref{d:rem:1}.

\begin{proposition}
  \label{d:prop:3}
  Totally quasi-convex metric compact Hausdorff spaces are totally convex.
\end{proposition}

\begin{proof}
  Let \( X \) be a totally quasi-convex metric compact Hausdorff space, and
  let \( x_1, x_2 \in X \) and \( u_1, u_2 \in [0,\infty] \) such that \(
  d(x_1,x_2) = u_1+u_2 \).  Hence, for all \( \varepsilon > 0 \), there exists
  \( z \in X \) such that \( d(x_1,z) \leq u_1+\varepsilon \) and \( d(z,x_2)
  \leq u_2+\varepsilon \).  For every \(\varepsilon>0\), we consider the set
  \begin{align*} A_\varepsilon & = \{ y\in X \mid d(x_1,y) \leq
  u_1+\varepsilon \:\text{and}\: d(y,x_2) \leq u_2+\varepsilon \}\\
  &=d(x_1,-)^{-1}([0,u_1+\varepsilon])\cap d(-,x_2)^{-1}([0,u_2+\varepsilon]),
  \end{align*} since the metric \(d \colon X\times X\to [0,\infty]\) is lower
  semicontinuous, \(A_{\varepsilon}\) is closed in \(X\).  By hypothesis,
  \(A_{\varepsilon}\neq\varnothing\) for each \( \varepsilon > 0 \), and since
  \( A_\varepsilon \subseteq A_\delta \) whenever \( \varepsilon \leq \delta
  \), the collection of these sets is a filter basis.  Thus, we may choose an
  ultrafilter \(\mathfrak{x} \) which contains these sets.

  Since \( X \) is compact Hausdorff, \( \mathfrak{x} \) converges to (a
  unique) \( z \in X \). Hence, for every \(\varepsilon>0\), we have that
  \(z\in A_{\varepsilon}\), that is, \( d(x_1,z) \leq u_1+\varepsilon \) and
  \( d(z,x_2) \leq u_2 + \varepsilon \).  Thus, we obtain \( d(x_1,z) \leq u_1
  \) and \( d(z,x_2) \leq u_2 \).
\end{proof}

\begin{remark}
  As it is shown in \cite{HR18}, a metric compact Hausdorff space can be
  equivalently described as a metric space \(X\), equipped with a compact
  Hausdorff topology, so that the ultrafilter convergence \(\alpha \colon
  \mathcal{U}(X)\to X\) is a metric map.  This description allows for a
  slightly shorter argument in the proof of Proposition~\ref{d:prop:3}, based
  on Proposition~\ref{d:prop:2}.  If \(d(x_1,x_2)=u_1+u_2\), then there is
  some ultrafilter \(\mathfrak{z}\) on \(X\) with \(d(\doo{x_1},
  \mathfrak{z})\leq u_1\) and \(d(\mathfrak{z},\doo{x_2})\leq u_2\). Therefore
  \(d(x_1,\alpha(\mathfrak{z}))\leq u_1\) and
  \(d(\alpha(\mathfrak{z}),x_2)\leq u_2\).
\end{remark}

The following result follows from standard arguments in metric spaces.

\begin{proposition}
  \label{d:prop:1}
  Let \(X\) be a compact Hausdorff space equipped with a metric \(d \colon
  X\times X\to [0,\infty]\). Then the following assertions are equivalent.
  \begin{enumerate}[(i)]
    \item
      The upper/lower balls of the metric space \(X\) are open in the topology of \(X\).
    \item
      The symmetric open balls of \(X\) are open in the topology of \(X\).
    \item
      The map \(d\colon X\times X\to [0,\infty]\) is continuous with respect
      to the Euclidean topology of \([0,\infty]\).
  \end{enumerate}
\end{proposition}

We denote by \(\MKH\) the category of metric compact Hausdorff spaces and
metric and continuous maps. Both the max-metric and the plus-metric endow
\(\MKH\) with a symmetric monoidal structure; however, as we show in the next
section, \(\MKH\) not closed for either of these structures.  For a metric
space \(X\), the function \(d \colon X\times X\to [0,\infty]\) is continuous
with respect to the topology on \(X\) induced by the metric and the Euclidean
topology on \([0,\infty]\), therefore every compact separated metric space is
a metric compact Hausdorff space when equipped with the topology induced by
the metric.  We denote the full subcategory of \(\Met\) formed by all compact
separated metric spaces by \(\CompMet\), and therefore we have a full
embedding
\begin{displaymath}
  \CompMet \longrightarrow \MKH.
\end{displaymath}

From Proposition~\ref{d:prop:1} and Remark~\ref{d:rem:2} we obtain the
following observation.

\begin{lemma}
  \label{d:lem:1}
  A separated metric compact Hausdorff space is compact metric if and only if
  any (and hence all) of the three assertions of Proposition~\ref{d:prop:1}
  hold.
\end{lemma}

We conclude this section with a metric version of Urysohn's Lemma essentially
obtained in \cite[Theorem~2.4]{Mat00} (see also \cite[Theorem~1.6]{BYaa13}).

\begin{theorem}
  \label{d:thm:1}
  Let \(X\) be a metric compact Hausdorff space and \(A,B\subseteq X\) closed
  subsets.  Put \(r\coloneq \inf\{d(x,y)\mid x\in A,y\in B\}\).  Then there
  exists a continuous metric map \(f \colon X \to [0,r]\) with \(A\subseteq
  f^{-1}(0)\) and \(B\subseteq f^{-1}(r)\).
\end{theorem}

\begin{corollary}
  \label{d:cor:2}
  The set \(\{[0,r] \mid 0<r<\infty\}\) and the space \([0,\infty]\) are
  initially dense in \(\MKH\).
\end{corollary}

\section{Main results}
\label{rui:sect:exp}
To treat various monoidal structures on \(\Met\) and \(\MKH\) simultaneously,
we consider here an arbitrary ordered commutative monoid structure \(\otimes\)
on \([0,\infty]\) with neutral element \(0\).  Furthermore, we assume that
\(\otimes\) satisfies \textit{Mulholland's inequality}\footnote{This condition
simply states that $\otimes \colon \Met \times \Met \to \Met$ is a lax
monoidal functor w.r.t to the tensor induced by the plus-metric. In
\cite{AFS06}, this means $+$ \textit{dominates} $\otimes$.} \cite{Mot36,
Mul49} (a generalization of Minkowski's inequality (Example~\ref{d:ex:3})):
\begin{displaymath}
  (u_1+ v_1)\otimes(u_2+ v_2)
  \leq
  (u_1\otimes u_2)+(v_1\otimes v_2),
\end{displaymath}
for all $u_1,u_2,v_1,v_2\in [0,\infty]$.

Indeed, this inequality ensures that the function \( d \colon (X\times Y)
\times (X\times Y) \to [0,\infty] \) defined by
\begin{displaymath}
  d\big((x,y),(x',y')\big)\coloneq X(x,x')\otimes Y(y,y'),
\end{displaymath}
for all \(x,x'\in X\) and \(y,y'\in Y\), is a metric on the product space \( X
\times Y \).  Naturally, the reflexivity condition holds, and so does the
triangle inequality: we have
\begin{align*}
  d\big((x,y),(x',y')\big)
    &= X(x,x') \otimes Y(y,y') \\
    &\leq \big( X(x,x'') + X(x'',x') \big)
           \otimes \big( Y(y,y'') + Y(y'',y') \big) \\
    &\leq \big( X(x,x'') \otimes Y(y,y'') \big)
           + \big( X(x'',x') \otimes Y(y'',y') \big) \\
    &= d\big((x,y),(x'',y'')\big) + d\big((x'',y''),(x',y')\big)
\end{align*}
for all \( x,x',x'' \in X \) and \( y,y',y'' \in Y \). We denote this metric
space by \( X \otimes Y \), and we note that this defines a symmetric monoidal
structre on the category \( \Met \).

\begin{lemma}
  If \( \otimes \) is an ordered commutative monoid structure on \( [0,\infty]
  \) satisfying Mulholland's inequality, then:
  \begin{itemize}
    \item
      we have \( a+b \geq a \otimes b \geq \max(a,b) \) for all \(a,b\in
      [0,\infty]\),
    \item
      \( \otimes \) defines a metric map \( [0,\infty] \otimes [0,\infty] \to
      [0,\infty] \).
  \end{itemize}
\end{lemma}
\begin{proof}
  For the inequalities, simply note that, for all \( a,b \in [0,\infty] \),
  \begin{equation*}
    a+b = (a \otimes 0) + (0 \otimes b) \geq (a + 0) \otimes (0+b) = a \otimes
    b,
  \end{equation*}
  and that monotonicity of \( \otimes \) guarantees \( a\otimes b \geq
  a\otimes 0 = a \), and \( a \otimes b \geq 0 \otimes b = b \).

  The second reduces to proving that \( \mu(a\otimes c,b \otimes d) \leq
  \mu(a,b) \otimes \mu(c,d) \)\footnote{Category theorists might recognise
  this as the \textit{mate} of Mulholland's inequality.} for all \( a,b,c,d
  \in [0,\infty] \). This follows by monotonicity and Mulholland's
  inequality.
\end{proof}

This yields numerous consequences. For one, we have that the identity
functions \( A \oplus B \to A \otimes B \to A \times B \) are metric maps, and
therefore the topology on \( A \otimes B \) induced by the metric is the
product topology. It also follows that \( \otimes \colon [0, \infty] \times
[0,\infty] \to [0,\infty] \) is a continuous map, and therefore this tensor
product defines a symmetric monoidal structure on \(\MKH\) as well.

Additionally, we also find that the map \(-\otimes u \colon [0,\infty]\to
[0,\infty]\) preserves infima for all \(u\in [0,\infty]\).  Thus, \(-\otimes
u\) has a right adjoint with respect to the ``greater or equal relation'' on
\([0,\infty]\), which we denote by \([u,-]_{\otimes}\colon [0,\infty]\to
[0,\infty]\). In other words, Mulholland's inequality ensures \(\otimes\) is a
commutative quantale on \( [0,\infty] \) with order relation \(\geq\) and unit
\(0\) (see \cite{Ros90a}).

We call a metric (compact Hausdorff) space \emph{\(\otimes\)-exponentiable}
whenever \(-\otimes X\) has a right adjoint, denoted by \([X,-]_{\otimes}\).

\begin{example}
  \label{d:ex:3}
  \begin{enumerate}[label=(\alph*)]
    \item
      We may consider \(\otimes=+\) and \(\otimes=\max\).
    \item
      \label{r:minkowski}
      For \(1<p<\infty\), we may consider \(u\otimes_p v=\sqrt[p]{u^p+v^p}\),
      thanks to \emph{Minkowski's inequality}:
      \begin{displaymath}
        \sqrt[p]{{(u_1 + v_1)}^p + {(u_2 + v_2)}^p}
        \leq
        \sqrt[p]{u_1^p + u_2^p} + \sqrt[p]{v_1^p + v_2^p},
      \end{displaymath}
      which additionally states that the equality holds if and only if \(
      u_1v_2 = u_2v_1 \).

      We point out that \( \otimes_1 = + \) and \( \lim_{p \to \infty}
      \otimes_p = \max \).
    \item
      From \cite[Theorem 1]{Mul49}, we deduce that if \( g \colon [0,\infty]
      \to [g(0),\infty] \) is a strictly increasing, continuous, convex and
      \textit{geoconvex}\footnote{A real-valued function $ g $ is geoconvex if
      $ g(\sqrt{xy}) \leq \sqrt{g(x)g(y)} $ for all $x,y$ in the domain.}
      function, then the binary operation defined by \( x \otimes y =
      g^{-1}\big(g(x)+g(y)-g(0)\big) \) is an ordered commutative monoid
      structure on \( [0,\infty] \) satisfying Mulholland's inequality.

      We have \( g(x) = e^x \) as an example that is not a power function.
    \item
      \cite[Table 2.6]{AFS06} lists several more examples\footnote{Not all;
      examples can only be found among the strict copulas.}, as one may check.
  \end{enumerate}
\end{example}

\begin{remark}
  By general adjunction rules, there is a bijection
  \([X,Y]_{\otimes}\cong\Met(X,Y)\), and the diagrams
  \begin{displaymath}
    \begin{tikzcd}
      {[X,Y]}_{\otimes}\otimes X\ar{d}[swap]{\cong}\ar{dr}{\text{counit}} \\
      \Met(X,Y)\otimes X
      \ar{r}[swap]{\ev} & Y,
    \end{tikzcd}
    \qquad
    \begin{tikzcd}
      & {[X,Y\otimes X]}_{\otimes}\ar{d}{\cong}\\
      Y\ar{r}[swap]{\eta}\ar{ur}{\text{unit}} & \Met(X,Y\otimes X)
    \end{tikzcd}
  \end{displaymath}
  commute. Here, the canonical map \(\eta \colon Y\to \Met(X,Y\otimes X)\)
  sends \(y\in Y\) to the metric map defined by \(x\mapsto (y,x)\). For an
  \(\otimes\)-exponentiable metric compact Hausdorff space, the same
  assertions hold, with \(\MKH(X,Y)\) in place of \(\Met(X,Y)\).  Henceforth,
  we take the underlying set of \([X,Y]_{\otimes}\) to be equal to
  \(\Met(X,Y)\) respectively \(\MKH(X,Y)\), and the unit and the counit to be
  given by \(\eta\) and \(\ev\), respectively.
\end{remark}

Below we give a slight generalisation of Theorem~\ref{d:thm:2}.

\begin{theorem}
  \label{d:thm:3}
  If \(X\) be a metric space, then the following are equivalent:
  \begin{enumerate}[(i)]
  \item \(X\) is \(\otimes\)-exponentiable in \(\Met\).
  \item\label{d:item:4} For all \( x_1, x_2 \in X \) and all \( u_1, u_2 \in [0,\infty]\),
    \begin{displaymath}
      \inf_{z\in X}\big((d(x_1,z)\otimes u_1)+(d(z,x_2)\otimes u_2)\big)
      \leq d(x_1,x_2)\otimes(u_1+u_2).
    \end{displaymath}
  \end{enumerate}
  If \(X\) is \(\otimes\)-exponentiable in \(\Met\), then the metric of
  \([X,Y]_{\otimes}\) is given by
  \begin{equation}
    \label{d:eq:1}
    d(f,g)=\sup_{x,x'\in X}[d(x,x'),d(f(x),g(x'))]_{\otimes},.
  \end{equation}
\end{theorem}
\begin{proof}
  For a \(\otimes\)-exponentiable metric space \( X \), its counit \(\ev \colon
  [X,Y]_{\otimes}\otimes X\to Y\) is a metric map, hence
  \begin{displaymath}
    d(f,g)\geq
    \underbrace{\sup_{x,x'\in X}[d(x,x'),d(f(x),g(x'))]_{\otimes}}_{\eqcolon u}.
  \end{displaymath}
  Consider now the metric space \(Z=\{0,1\}\) with \(d(0,1)=u\) and \(d(1,0)=\infty\).
  Then the map
  \begin{displaymath}
    X\otimes Z \longrightarrow Y,\quad
    (x,0)  \longmapsto f(x),\;
    (x,1)  \longmapsto g(x)
  \end{displaymath}
  is metric, hence also the map
  \begin{displaymath}
    Z \longrightarrow [X,Y]_{\otimes},\quad
    0 \longmapsto f,\;
    1 \longmapsto g
  \end{displaymath}
  is metric, and therefore \(u\geq d(f,g)\).  We conclude that a metric space
  \(X\) is \(\otimes\)-exponentiable in \(\Met\) if and only if the distance
  defined by \eqref{d:eq:1} satisfies the triangle inequality, for all metric
  spaces \(Y\).  With the same proof as in \cite{CH06} one verifies that this
  is indeed equivalent to the assertion~\ref{d:item:4} above.  In fact, an
  easy calculation shows that \ref{d:item:4} guarantees the \(d\) satisfies
  the triangle inequality, for all metric spaces \(Y\).

  Assume now that \(X\) is \(\otimes\)-exponentiable. We consider the maps
  \begin{align*}
    \varphi_1 \colon X & \longrightarrow [0,\infty],
    \; x \longmapsto d(x_1,x),\\
    \varphi_2 \colon X & \longrightarrow [0,\infty],
    \; x \longmapsto d(x_1,x)\otimes u_1,\\
    \varphi_3 \colon X & \longrightarrow [0,\infty],
    \; x \longmapsto \inf_{z\in X}\big((d(x_1,z)\otimes u_1)+(d(z,x)\otimes u_2)\big).
  \end{align*}
  for each \( x_1, x_2 \in X \) and \( u_1, u_2 \in [0,\infty]\).  Then
  \(\varphi_1,\varphi_2,\varphi_3\) are in fact metric maps since, for all
  \(x\in X\), all \(u\in [0,\infty]\) and all sets \(I\), the maps
  \begin{displaymath}
    d(x,-)\colon X\longrightarrow [0,\infty],\quad
    -\otimes u \colon [0,\infty] \longrightarrow [0,\infty],\quad
    \inf \colon [0,\infty]^I \longrightarrow [0,\infty]
  \end{displaymath}
  are metric.  For all metric maps \(\varphi,\varphi' \colon X\to [0,\infty]\)
  and \(u\in [0,\infty]\), we have \(u\geq d(\varphi,\varphi')\) if and only
  if, for all \(x,x'\in X\), \(u\otimes d(x,x') \geq \mu(\varphi(x),
  \varphi'(x'))\), which in turn is equivalent to \((u\otimes
  d(x,x´))+\varphi(x)\geq \varphi'(x')\), for all \(x,x'\in X\).  From that we
  conclude that \(u_1\geq d(\varphi_1,\varphi_2)\) and \(u_2\geq
  d(\varphi_2,\varphi_3)\).  Since \(d\) satisfies the triangle inequality,
  \(u_1+u_2\geq d(\varphi_1,\varphi_3)\).  Therefore
  \begin{align*}
    (u_1+u_2)\otimes d(x_1,x_2)
    &=((u_1+u_2)\otimes d(x_1,x_2))+\varphi_1(x_1)\\
    &\geq\varphi_3(x_2)
    =\inf_{z\in X}\big((d(x_1,z)\otimes u_1)+(d(z,x_2)\otimes u_2)\big).
    \qedhere
  \end{align*}
\end{proof}

We remark that it is even possible to characterise \( \otimes
\)-exponentiability for a space \( X \) in terms of a more useful metric
condition on \( X \). Intuitively, \( X \) must simultaneously ``almost''
satisfy the converses of the triangle inequality and Mulholland's inequality.

To make this intuition precise, we fix the following notation for a subset \(
S \subseteq X \):
\begin{align*}
   B_\varepsilon(S)
    &= \{ x \in X \mid \exists s \in S\:: d(x,s) < \varepsilon \text{
     and } d(s,x) < \varepsilon \}, \\
   \overline{B}_\varepsilon(S)
    &= \{ x \in X \mid \exists s \in S\:: d(x,s) \leq \varepsilon \text{
     and } d(s,x) \leq \varepsilon \}.
\end{align*}
Additionally, we consider the continuous function
\begin{displaymath}
  \Omega_{d,u_1,u_2} \colon {[0,d]} \longrightarrow {[0,\infty]},\;
  r \longmapsto r \otimes u_1 + (d-r) \otimes u_2
\end{displaymath}
and note that we have \( \Omega(r) \geq d \otimes (u_1+u_2) \) for all \(r\),
by Mulholland's inequality. With this in mind, we define
\begin{equation*}
  M_{d,u_1,u_2} = \Omega^{-1}\big(d \otimes (u_1+u_2)\big)
                = \{ r \in [0,d] \mid r \otimes u_1 + (d-r) \otimes u_2
                      = d \otimes (u_1 + u_2) \}.
\end{equation*}
Note that \( M_{0,u_1,u_2} = \{0\} \), \( M_{\infty,u_1,u_2} = [0,\infty]
\), \( M_{d,0,0} = [0,d] \), and \(  M_{d,\infty,u_2} = M_{d,u_1,\infty} =
[0,d] \). These are the cases \( d = 0, d=\infty \), \( u_1+u_2=0 \) and \(
u_1+u_2 = \infty \), respectively. We call these the \textit{boundary
conditions}.

\begin{theorem}
  \label{thm:blob-convexity}
  Assume that \( a \otimes x < b \otimes x \) whenever \( a < b \) and \( x <
  \infty \).\footnote{Such binary operations are called \textit{strict} in
  \cite{AFS06}.} The following are equivalent for a metric space \( X \):
  \begin{enumerate}[(\roman*)]
    \item
      \label{enum:exponentiable}
      \( X \) is \( \otimes \)-exponentiable.
    \item
      \label{enum:blob-convexity}
      for all \( x, y \in X \), all \( u_1,u_2 \in [0,\infty] \) and all \(
      \varepsilon > 0 \), there exists \( z \in X \) with \( d(x,z) + d(z,y) \leq
      d(x,y) + \varepsilon \) and \( d(x,z) \in \overline{B}_\varepsilon(M_{d(x,y),
      u_1, u_2}) \).
  \end{enumerate}
\end{theorem}

\begin{proof}
  We use Theorem~\ref{d:thm:3}.  To prove \ref{enum:exponentiable} implies
  \ref{enum:blob-convexity}, let \( x,y \in X \), let \( u_1, u_2 \in
  [0,\infty] \), and let \( \varepsilon > 0 \).  Under the boundary conditions,
  there is nothing to show; hence, we assume \( 0 < d(x,y) < \infty \), and \(
  0 < u_1 + u_2 < \infty \).

  We write \( M = M_{d(x,y),u_1,u_2} \); if \( B_{\varepsilon}(M) = [0,d(x,y)] \),
  there is nothing to check; else, we take
  \begin{equation*}
    \eta_1 = \min \Omega|_{[0,d(x,y)]\setminus B_\varepsilon(M)},
    \qquad \eta_2 = \big(d(x,y) + \varepsilon\big) \otimes (u_1+u_2),
  \end{equation*}
  highlighting that if \( \Omega(a) \leq \eta_1 \), then \( a \in
  \overline{B}_\varepsilon(M) \).

  Taking \( \omega = \min(\eta_1,\eta_2) - d(x,y)\otimes(u_1+u_2) > 0 \), by
  \( \otimes \)-exponentiability we find \( z \in X \) with
  \begin{equation*}
    d(x,z) \otimes u_1 + d(z,y) \otimes u_2 \leq d(x,y) \otimes (u_1+u_2) +
    \omega = \min(\eta_1,\eta_2).
  \end{equation*}
  We have \( d(x,y) - d(x,z) \leq d(z,y) \), from which we deduce that \(
  \Omega(d(x,z)) \leq \eta_1 \), and hence \( d(x,z) \in
  \overline{B}_\varepsilon(M) \).
  Moreover, we also have
  \begin{equation*}
    \big(d(x,z) + d(z,y)\big) \otimes (u_1+u_2)
      \leq d(x,z) \otimes u_1 + d(z,y) \otimes u_2 \leq \eta_2,
  \end{equation*}
  from which we deduce that \( d(x,z) + d(z,y) \leq d(x,y) + \varepsilon \).

  To show \ref{enum:blob-convexity} implies \ref{enum:exponentiable}, let \(
  x,y \in X \), \( u_1,u_2 \in [0,\infty] \), and let \( \eta > 0 \). With
  nothing to prove under the boundary conditions, we assume \( 0 < d(x,y) <
  \infty \), and \( 0 < u_1 + u_2 < \infty \).

  Let \( \varepsilon_1 = \frac{\eta}{2} > 0 \).  Since \( \Omega \) is a
  continuous function defined on a compact interval, there exists \(
  \varepsilon_2 > 0 \) such that if \( r \in \overline B_{\varepsilon_2}(M) \), then
  \( \Omega(r) < d(x,y) \otimes (u_1+u_2) + \varepsilon_1 \). Indeed, for all $x
  \in M$, there exists $ \delta > 0 $ such that for all $y$ with $ |x-y| <
  \delta $, we have $\Omega(y) - \Omega(x) < \varepsilon_1 $. Now, by
  compactness, we find finitely many $\delta$; take \( \varepsilon_2 \) to be
  the smallest.

  So, with \( \varepsilon = \min(\varepsilon_1,\varepsilon_2) > 0\), via 
  \ref{enum:blob-convexity} we find \( z \in X \) such that \( d(x,z) \in
  \overline{B}_{\varepsilon_2}(M) \) and \( d(x,z)+d(z,y) \leq d(x,y) +
  \varepsilon_1 \).

  Thus, we have that
  \begin{align*}
    d(x,z) \otimes u_1 + d(z,y) \otimes u_2
      &\leq d(x,z) \otimes u_1
              + \big(d(x,y)-d(x,z)+\tfrac{\eta}{2} \big) \otimes u_2 \\
      &\leq \Omega\big(d(x,z)\big) + \tfrac{\eta}{2} \leq d(x,y) \otimes (u_1+u_2) + \eta,
  \end{align*}
  as desired.
\end{proof}

As a corollary, we characterise \( \otimes_p \)-exponentiability:

\begin{corollary}
  \label{cor:minkowski}
  If \( p > 1 \), a metric space \(X \) is \( \otimes_p \)-exponentiable if
  and only if \( X \) is totally quasi-convex.
\end{corollary}
\begin{proof}
  As noted in Example~\ref{d:ex:3}\ref{r:minkowski}, we have \( M_{d,u_1,u_2}
  = \{ d \cdot \frac{u_1}{u_1+u_2} \} \) for \( 0 < u_1 + u_2 < \infty \) and
  \( 0 < d < \infty \), by Minkowski's inequality.  One then immediately
  checks that Theorem~\ref{thm:blob-convexity}\ref{enum:blob-convexity} is
  equivalent to total quasi-convexity.
\end{proof}

\begin{corollary}
  Under the hypothesis of Theorem~\ref{thm:blob-convexity}, if
  \(M_{d,u_1,u_2}=\varnothing\) outside the boundary conditions, then a metric
  space is \(\otimes\)-exponentiable if and only if it is discrete.
\end{corollary}

\begin{remark}
  We conjecture that \(M_{d,u_1,u_2}=\varnothing\) outside the boundary
  conditions for all tensors not of the form \(\otimes_p\).
\end{remark}

In the remainder of this section we analyse \(\otimes\)-exponentiability for metric
compact Hausdorff spaces; we start by providing a necessary condition.

\begin{proposition}
  \label{d:prop:4}
  Let \(X\) be a metric compact Hausdorff space.
  Then (\ref{d:item:1})\(\implies\)(\ref{d:item:2})\(\implies\)(\ref{d:item:3}) below.
  \begin{enumerate}
    \item
      \label{d:item:1}
      \(X\) is \( \otimes \)-exponentiable on \(\MKH\).
    \item
      \label{d:item:2}
      For every metric compact Hausdorff space \(Y\), \(\MKH(X,Y)\) is a
      closed subspace of the product space \(\prod_{x\in X}Y\).
    \item
      \label{d:item:3}
      The lower/upper open balls are open in the topology of \(X\).
  \end{enumerate}
\end{proposition}

\begin{proof}
  Assume first (\ref{d:item:1}). For every \(x\in X\), \(\ev_x \colon
  [X,Y]_{\otimes}\to Y\) is an arrow in \(\MKH\), in particular it is
  continuous, hence we obtain a canonical injective and continuous map
  \(\MKH(X,Y)\to \prod_{x\in X}Y\). By compactness, this map is a closed
  embedding.

  Assume now (\ref{d:item:2}).  Hence, for each \(x\in X\), the map
  \(\pi_{x}\colon \MKH(X,[0,\infty])\to [0,\infty],\, f\mapsto f(x)\) is
  continuous.  For every \(\varepsilon<\varepsilon_0\) and \(U\subseteq X\)
  open with \(\BRight{x_0}{\varepsilon_0}\subseteq U\) we put
  \begin{displaymath}
    \mathcal{C}_{U,\varepsilon} \coloneqq
      \{h \colon X\to [0,\infty] \text{ in }\MKH 
        \mid h(x_0)=0\text{ and } h(z)\geq\varepsilon
              \text{ for all }z\in U^{\complement}\}.
  \end{displaymath}
  By Theorem~\ref{d:thm:1}, \(\mathcal{C}_{U,\varepsilon}\neq\varnothing\);
  moreover,
  \begin{displaymath}
    \mathcal{C}_{U,\varepsilon}\cap \mathcal{C}_{V,\delta}\supseteq
    \mathcal{C}_{U\cap V,\max\{\varepsilon,\delta\}}.
  \end{displaymath}
  Therefore there exists an ultrafilter \(\mathfrak{x}\) on
  \(\MKH(X,[0,\infty])\) such that \(\mathcal{C}_{U,\varepsilon}\in
  \mathfrak{x}\) for all \(\varepsilon<\varepsilon_0\) and \(U\subseteq X\)
  open with \(\BRight{x_0}{\varepsilon_0}\subseteq U\). By compactness,
  \(\mathfrak{x}\) converges to some \(f\in\MKH(X,[0,\infty])\). We have the
  following facts.
  \begin{enumerate}
    \item 
      \(f(x_0)=0\). In fact, for every \(\mathcal{A}\in \mathfrak{x}\),
      \(\mathcal{A}\cap \mathcal{C}_{X,\frac{\varepsilon_0}{2}}\neq
      \varnothing\), that is, there is some \(k\in \mathcal{A}\) with
      \(k(x_0)=0\). Therefore \(\pi_{x_0}(\mathfrak{x})\) is the principal
      ultrafilter generated by \(0\). Since \(\pi_{x_0}(\mathfrak{x})\to
      f(x_0)\), we conclude that \(f(x_0)=0\).
    \item 
      \(f(x)<\varepsilon_0\) for every \(x \in
      \BRight{x_{0}}{\varepsilon_{0}}\). In fact, since \(f\) is a metric map,
      \begin{displaymath}
        f(x)=f(x)-f(x_0)\leq X(x_0,x)<\varepsilon_0.
      \end{displaymath}
    \item  
      \(f(y)\geq\varepsilon_0\) for every
      \(y\notin\BRight{x_{0}}{\varepsilon_{0}}\). To see this, note first
      that, by Theorem~\ref{d:thm:1}, for every \(\varepsilon<\varepsilon_0\)
      there exists an arrow \(h_{\varepsilon}\colon X\to [0,\infty]\) in
      \(\MKH\) with \(h_{\varepsilon}(x_0)=0\) and \(h_{\varepsilon}(y)\geq
      \varepsilon\).  Put
      \begin{displaymath}
        U=\bigcup_{\varepsilon<\varepsilon_0}h_{\varepsilon}^{-1}([0,\varepsilon[).
      \end{displaymath}
      Then \(U\subseteq X\) is open,
      \(\BRight{x_{0}}{\varepsilon_{0}}\subseteq U\) and \(y\notin U\).

      Let now \(\varepsilon<\varepsilon_0\). For every \(\mathcal{A}\in
      \mathfrak{x}\), \(\mathcal{A}\cap
      \mathcal{C}_{U,\varepsilon}\neq\varnothing\), that is, there is some
      \(k\in \mathcal{A}\) with \(k(y)\geq \varepsilon\). Therefore
      \([\varepsilon,\infty]\in \pi_y(\mathfrak{x})\), since
      \(\pi_y(\mathfrak{x})\to f(y)\), we conclude that \(f(y)\geq
      \varepsilon\). Consequently, \(f(y)\geq \varepsilon_0\).
  \end{enumerate}
  Putting everything together, we have shown that
  \(\BRight{x_{0}}{\varepsilon_{0}}=f^{-1}([0,\varepsilon_0[)\) is open in
  \(X\).  A similar argument shows that \(\BLeft{x_0}{\varepsilon_0}\) is open
  in \(X\).
\end{proof}

Using Lemma~\ref{d:lem:1} and Proposition~\ref{d:prop:4}, we deduce
immediately the following result.

\begin{theorem}
  Let \(X\) be a \(\otimes\)-exponentiable separated metric compact Hausdorff space.
  Then \(X\) is a compact metric space.
\end{theorem}

Using Example~\ref{d:ex:2}, we conclude the following fact for ordered compact
Hausdorff spaces.

\begin{corollary}
  \label{d:cor:1}
  Let \(X\) be an exponentiable ordered compact Hausdorff space. Then, for
  every \(x\in X\), the sets \(\downc x=\{y\in X\mid y\leq x\}\) and \(\upc
  x=\{y\in X\mid x\leq y\}\) are open in \(X\).  In particular, if \(X\) is
  antisymmetric, then \(X\) is finite.
\end{corollary}

Now that we have a better understanding of the topological structure of
\(\otimes\)-exponentiable (separated) metric compact Hausdorff spaces, we
expose their metric structure.

\begin{proposition}
  A \(\otimes\)-exponentiable metric compact Hausdorff space is also
  \(\otimes\)-exponentiable in \(\Met\).
\end{proposition}
\begin{proof}
  The underlying set of \([X,Y]_{\otimes}\) is necessarily given by
  \begin{math}
    \MKH(X,Y),
  \end{math}
  and the same argument as in the proof of Theorem~\ref{d:thm:3} shows that
  the distance of \([X,Y]_{\otimes}\) is given by~\eqref{d:eq:1}, using that
  the space \(Z\) is finite.  In fact, as in the proof of
  Theorem~\ref{d:thm:3}, we show that an exponentiable metric compact space
  satisfies Condition~\ref{d:item:4} of Theorem~\ref{d:thm:3}.  First observe
  that, since \(X\) is \(\otimes\)-exponentiable, the metric \(d \colon
  X\times X\to [0,\infty]\) is even continuous (see
  Propositions~\ref{d:prop:4} and \ref{d:prop:1}).  Since \(\otimes \colon
  [0,\infty]\times [0,\infty]\to [0,\infty]\) is continuous, we obtain
  immediately that the metric functions \(\varphi_1,\varphi_2\colon X\to
  [0,\infty]\) defined in the proof of Theorem~\ref{d:thm:3} are continuous.
  Regarding \(\varphi_3\colon X\to [0,\infty]\), observe first that the map
  \begin{displaymath}
    \psi \colon X\times X \longrightarrow [0,\infty],\;
    (z,x) \longmapsto (d(x_1,z)\otimes u_1)+(d(z,x)\otimes u_2)
  \end{displaymath}
  is continuous.  With \([0,\infty]^X\) denoting the exponential in the
  category of topological spaces and continuous maps (which can be taken since
  \(X\) is compact Hausdorff), its exponential mate
  \begin{displaymath}
    \overline{\psi}\colon X \longrightarrow [0,\infty]^X
  \end{displaymath}
  is also continuous.  Furthermore, \(\inf \colon [0,\infty]^X\to [0,\infty]\)
  is continuous since \(X\) is compact  (see the proof of
  \cite[Proposition~4.1]{AH25}), therefore
  \(\varphi_3=\inf\cdot\overline{\psi}\) is continuous.  As in the proof of
  Theorem~\ref{d:thm:3}, we conclude that
  \begin{displaymath}
    (u_1+u_2)\otimes d(x_1,x_2)\geq
    \inf_{z\in X}\big((d(x_1,z)\otimes u_1)+(d(z,x_2)\otimes u_2)\big)
  \end{displaymath}
  for all  \( x_1, x_2 \in X \) and all \( u_1, u_2 \in [0,\infty]\).
\end{proof}

We will now also see that the assertion~\ref{d:item:4} of
Theorem~\ref{d:thm:3} implies \(\otimes\)-exponentiabilty in \(\MKH\), for
compact separated metric spaces.

\begin{lemma}
  Let \(X\) be a metric space and \(Y\) be a metric compact Hausdorff space.
  Then the subset \(\Met(X,Y)\subseteq \prod_{x\in X}Y\) is closed in the product
  space \(\prod_{x\in X}Y\).
\end{lemma}
\begin{proof}
  Let \(f\in Y^{|X|}\) with \(f\notin\Met(X,Y)\). Then there exist
  \(x_1,x_2\in X\) with \(d(x_1,x_2)<d(f(x_1),f(x_2))\). Put
  \(\varepsilon=d(x_1,x_2)\), then \(d(f(x_1),f(x_2))=\varepsilon+\eta\) with
  \(\eta>0\).  By Theorem~\ref{d:thm:1}, there exists \(h \colon
  X\to[0,\infty]\) in \(\MKH\) with \(h(f(x_1))=0\) and
  \(h(f(x_2))=\varepsilon+\eta\).  Put \(U=h^{-1}([0,\frac{\eta}{2}[)\) and
  \(V=h^{-1}(]\varepsilon+\frac{\eta}{2},\infty])\), and note that \(U\) and
  \(V\) are open subsets of \(X\) with \(f(x_1)\in U\) and \(f(x_2)\in V\).
  Moreover, for all \(y_1\in U\) and \(y_2\in V\), \(d(y_1,y_2)\geq
  h(y_2)-h(y_1)>\varepsilon\).  Therefore
  \(W=\pi_{x_1}^{-1}(U)\cap\pi_{x_2}^{-1}(V)\) is an open neighbourhood of
  \(f\) with \(W\cap\Met(X,Y)=\varnothing\).  We conclude that \(\Met(X,Y)\)
  is closed in \(Y^{|X|}\).
\end{proof}

\begin{corollary}
  For every metric space \(X\) and every metric compact Hausdorff space \(Y\),
  we have that \(\Met(X,Y)\) is a metric compact Hausdorff space with the
  structure inherited from \(\prod_{x\in X}Y\).
\end{corollary}

\begin{lemma}
  If \(X\) is a compact separated metric space and \(Y\) is a metric compact
  Hausdorff space, then the canonical map \(\eta \colon Y \to \Met(X,Y\otimes
  X)\) is continuous.
\end{lemma}

\begin{lemma}
  \label{d:lem:2}
  Let \(X\) and \(Y\) be compact separated metric spaces. Then the map
  \begin{math}
    \ev \colon \Met(X,Y)\otimes X \to Y
  \end{math}
  is continuous.
\end{lemma}
\begin{proof}
  Let \(f \colon X\to Y\) be a metric map, \(x\in X\) and \(V\subseteq Y\) an
  open neighbourhood of \(f(x)\). Since \(Y\) is a compact metric space, there
  is some \(\varepsilon>0\) so that \(\BBoth{f(x)}{2\varepsilon}\subseteq V\).
  Then \(W\coloneqq\pi_x^{-1}(\BBoth{f(x)}{\varepsilon})\times
  \BBoth{x}{\varepsilon}\) is an open neighbourhood of \((f,x)\) and, for all
  \((h,z)\in W\),
  \begin{displaymath}
    d(h(z),f(x))\leq d(h(z),h(x))+d(h(x),f(x))\leq 2\varepsilon
  \end{displaymath}
  and \(d(f(x),h(z))\leq 2\varepsilon\).  Therefore \(\ev(W)\subseteq V\).
\end{proof}

With the help of Corollary~\ref{d:cor:2}, we extend Lemma~\ref{d:lem:2} to all
metric compact Hausdorff spaces \(Y\).

\begin{lemma}
  Let \(X\) be a compact separated metric space and \(Y\) be a metric compact
  Hausdorff space.  Then the map
  \begin{math}
    \ev \colon \Met(X,Y)\otimes X \to Y
  \end{math}
  is continuous.
\end{lemma}
\begin{proof}
  Let \(\varphi \colon Y\to [0,r]\) be a continuous and metric map.  First
  observe that map
  \begin{displaymath}
    \Met(X,\varphi)\otimes 1_X \colon \Met(X,Y)\otimes X
    \longrightarrow
    \Met(X,[0,r])\otimes X,\quad
    (f,x) \longmapsto (\varphi f,x)
  \end{displaymath}
  is continuous (and metric). Secondly, the diagram
  \begin{displaymath}
    \begin{tikzcd} 
      \Met(X,Y)\otimes X 
      \ar{r}{\ev} 
      \ar{d}[swap]{\Met(X,\varphi)\otimes 1_{X}} 
      & Y 
      \ar{d}{\varphi} \\
      \Met(X,[0,r])\otimes X 
      \ar{r}[swap]{\ev} 
      & {[0,r]} 
    \end{tikzcd}
  \end{displaymath}
  commutes.  Since \([0,r]\) is a compact separated metric space, \(\ev \colon
  \Met(X,[0,r])\otimes X \to [0,r]\) is continuous (Lemma~\ref{d:lem:2}) and
  therefore also \(\ev\cdot(\Met(X,\varphi)\otimes 1_X)\). Since the cone of
  all metric and continuous maps \(Y\to [0,r]\) separates the points of \(Y\)
  (Corollary~\ref{d:cor:2}), \(\ev \colon \Met(X,Y)\otimes X \to Y\) is also
  continuous.
\end{proof}

\begin{lemma}
  If \(X\) is a compact metric space, \(\otimes\)-exponentiable in \(\Met\),
  and \(Y\) is a metric compact Hausdorff space, then \([X,Y]_{\otimes}\) is a
  metric compact Hausdorff space.
\end{lemma}
\begin{proof}
  First note that, for all \(x,x'\in X\), the composite
  \begin{displaymath}
    [X,Y]_{\otimes}\times [X,Y]_{\otimes}
    \xrightarrow{\,\ev_x\times\ev_{x'}\,} Y\times Y \xrightarrow{\;d\;}
    [0,\infty] \xrightarrow{\,[X(x,x'),-]_{\otimes}\,} [0,\infty]
  \end{displaymath}
  is lower semicontinuous, and hence
  \begin{displaymath}
    \begin{tikzcd}
      {[X,Y]}_{\otimes}\times [X,Y]_{\otimes}
      \ar{r}\ar[bend left=15]{rr}{d}
      & \prod_{x,x'\in X}[0,\infty]
      \ar{r}[swap]{\sup}
      & {[0,\infty]}
    \end{tikzcd}
  \end{displaymath}
  is lower semicontinuous as well.
\end{proof}

Putting everything together:

\begin{theorem}
  \label{d:thm:4}
  We consider a continuous commutative quantale structure \(\otimes\) on
  \([0,\infty]\) with neutral element \(0\) satisfying Mulholland's
  inequality, that is,
  \begin{displaymath}
    (u_1+ u_2)\otimes(v_1+ v_2) \leq
    (u_1\otimes v_1)+(u_2\otimes v_2),\quad
    \text{for all $u_1,u_2,v_1,v_2\in [0,\infty]$.}
  \end{displaymath}
  Then, a separated metric compact Hausdorff space is
  \(\otimes\)-exponentiable in \(\MKH\) if and only if it is compact metric
  and \(\otimes\)-exponentiable in \(\Met\).
\end{theorem}

\begin{corollary}
  In particular, a separated metric compact Hausdorff space is
  \begin{itemize}
    \item 
      \(+\)-exponentiable in \(\MKH\) if and only if it is compact metric.
    \item 
      \( \otimes_p \)-exponentiable in \(\MKH\) if and only if it is compact
      metric and totally convex, for \( 1 < p \leq \infty \).
  \end{itemize}
\end{corollary}

\bibliographystyle{amsplain}

\end{document}